\documentstyle[11pt]{article}
\begin{document}
\title{\Large\bf  A new Binary Number Code and a Multiplier, \\
		based on 3 as semi-primitive root of 1 mod $2^k$ }
\author{N.F.Benschop, ~~Geldrop (NL), ~~3-aug-1997}
\date{----------- Patent US-5923888 (13jul99)-------------}
\maketitle

\section{Prior art}

The usual parallel array multipliers [1, p164] are much too powerful for their
purpose, to be shown as follows. Assuming without loss of generality a square array,
the known parallel $n \times n$ bit array multipliers all have a structure consisting
of two main parts. An {\bf input part} with a 2-dimensional array of $n$(hrz) + $n$(vrt)
bitlines, for the two $n$-bit input operands $x$ and $y$, with an AND-gate at each of
the $n^2$ bitline crossings (details for signed TC code are neglected here). And a
{\bf processing part}, which accumulates this pattern of $n^2$ bits to the required
$2n$-bit result, using an array of some $n^2$ Full-Adders (FA). Various types of
Adder-array exist, like a normal array of $n$ rows of $n$ FA's each (for a compact
layout and small silicon area), or the known 'Wallace tree' [1, p167]  (with an irregular
and larger layout but less delay), or anything between these extremes, trading-off total
delay and silicon area.
\\[1ex]
The inefficiency of the usual adder array hardware is easily seen as follows.
The adder array can add {\bf any} $n \times n$ bit pattern of $n^2$ bits (there
are $2^{n.n}$ patterns), while for multiplication of two $n$-bit operands {\bf only}
$2^{2n}$ of these are ever input and processed (each $n$-bit row or column is either
all 0's or a copy of one operand). So the hardware is used for processing only a
very small fraction $2^{n+n}/2^{n.n}$ of all possible input patterns it {\it could}
process. Clearly, the hardware is much too powerfull for its purpose, and is used
very inefficiently. Some recoding schemes have been applied in the past to improve
the efficiency of multipliers.
\\[1ex]
For instance in the known Booth multiplier [1, p198], each successive bit-pair of
one input operand has value range \{0, 1, 2, 3\}, where 3 is recoded as $-1+4$.
The -1 causes a subtraction of the other operand, while '+4', as positive carry into
the next bit-pair position, implies an addition there. The result is an effective
reduction of the logic depth in the add/subtract array, and a corresponding speed-up,
at the cost of a more complex recoding of one operand, and extra subtract hardware.
\\ A similar recoding scheme, but now for both operands, and based on a deeper
algebraic property of the powers of 3 in the semigroup of binary multiplication
$M$(.) mod $2^k$, will next be proposed.

\section{Proposed new binary number code}

A better structure might be found by using the algebraic properties of the closed
system (semigroup) of binary multiplication mod $2^k$, such as associativity
$a(bc)=(ab)c$, commutativity $ab=ba$, and the iterative sub-structures or {\it
iteration class} $a^*=\{a^i\}$ of all powers of any number $a$. Especially $a=3$,
which generates the maximum possible iteration class of order $2^{k-2}$, to be proven
next. Exploiting this {\it 3* property} makes multipliers much more efficient.
\\[1ex]
For $k \geq 3$ bits the powers of 3 generate half of the odd residues. In other
words, in binary coded residues: 3 is a semi-primitive root of unity. A new binary
number code based on this property simplifies binary multiplication, and in fact
translates it to addition, using base 3 logarithm for odd residues. The proof is best
given by first considering residues mod $p^k$ for prime $p>$2, and then taking $p$=2
as special case. Denote a cyclic group of order $n$ by $C_n$ or $C[n]$.
\\[1ex]
{\bf Lemma}:
~For prime $p>$2, the cyclic subgroup $B=(p+1)^*$ mod $p^k$ has order $p^{k-1}$.
\\[1ex]
{\bf Proof}: ~The group of $units ~G$ of all $n$ with $\{n^i$=1\} mod $p^k$ for some
$i>$0, is known to be cyclic. Its order $(p-1).p^{k-1}$ has two relative prime factors,
so $G=A \times B$ is a direct product of two cycles. Here $B=(p+1)^*$ because $(p+1)^p=p^2$+1
mod $p^3$, and by induction $(p+1)^{p^m}=p^{m+1}$+1 mod $p^{m+2}$. The period of $p$+1, the
smallest $x$ with $(p+1)^x=1$ mod $p^k$, implies $m+1=k$, so $m=k-1$, yielding period
$p^{k-1}$. No smaller $x$ yields 1 mod $p^k$ since $|B|$ has only divisors $p^s$.
~$\spadesuit$
\\[1ex]
{\bf Corollary}~( {\it binary 3* property} ):
~For $p$=2 we have $p$+1=3, and it is readily verified that 3 does not generate $-1$
mod $2^k$ if $k \geq 3$, since $(2+1)^2>2^3$ ~(in binary code $3^2$=1001), while
$(p+1)^2=p^2+2p+1<p^3$ for all $p>$2. The carry in binary code is the cause of this
phenomenon. In fact $B=C_2.C[2^{k-2}]$ is not cyclic, with {\it sign} 2-cycle $C_2=
\{-1,1\}$.
Then $|3^{*}|=2^{k-2}$,  with 3 generating only half of the odd numbers mod $2^k$;
the other half are their complements. So each non-zero residue is $n=\pm 3^i.2^j$
mod $2^k$, with $i<2^{k-2}$ and $j<k$, while $n=0$ for $j=k$. ~$\spadesuit$

\subsection{ Example }

For instance mod 32 ($k$=5) the cycle 3* =\{$ 3, 9, -5, -15, -13, -7, 11, 1 $\}
has period 8, while the remaining 8 odd numbers are their complements, with a
two-component decomposition $G=C_2.C_8$= \{$-1,1$\} x 3* for all 16 odd numbers,
which allows component-wise multiplication. The 5-bit binary codes of $3^i$
are shown in the next table, as well as for $p>2$ the lower significant digits of
$(p+1)^{p^m}$ in $p$-ary code. The logic structure of the few least
significant bits of $3^i$ is rather simple, as boolean functions of the $k-2$
exponent bits, but the higher order bits quickly increase in complexity, showing
no obvious structure.
\\[1ex]
{\bf Table 1}:
~The powers of 3 in binary code mod $2^5$, ~~~~~~~~~~and $(p+1)^{p^m}$ in p-ary code:
\begin{verbatim}
 i   3^i (bin)  3^i (dec)                              | p>2  (p+1)^i   i
1.    00011      3          Notice:  3^even = 1 mod 8  |-------------+----
2.    01001      9                   3^ odd = 3 mod 8  |        11     1
3.    11011     27 =  -5    so two bits  are fixed:    |     ..101     p
4.    10001     17 = -15     bit(2^0)= 1, bit(2^2)= 0  |  ....1001     p^2
5.    10011     19 = -13     hence:    |3*| = 2^k / 4  |.....10001     p^3
6.    11001     25 =  -7                               |
7.    01011     11
8.    00001      1

Table 2: ------------- Multiplier structure ------------
   Operands     a = sign(a) 3^i.2^j
                b = sign(b) 3^r.2^s         | sign(p)= XOR(signs)
   Product p= a.b = sign(p) 3^t.2^u  where: | t= i+r   mod 2^{k-2}
                                            | u= j+s   < k (saturate at k)
                                                            'overflow'
\end{verbatim}

\section{Application to multipliers}

By the corollary each residue is $n=\pm ~3^i.~2^j$ mod $2^k ~(k>$2) for a unique pair
($i,j$) of exponents, with  $0 \leq i<2^{k-2} ~~(k$-2 bits mantissa) and $0 \leq j<k$
(binlog $k$ bits), with $n$=0 iff $j$=$k$. ~This {\it 2.3-star} number code {\bf reduces
multiplication to addition} of exponent-pairs, because: $(3^i.2^j).(3^r.2^s)=
3^{i+r}.2^{j+s}$, and the 1-bit signs add (mod 2). The multiplier structure is summarized
in table 2: the product sign is the XOR or the operand signs, the exponents of 3 add mod
$2^{k-2}$ using only the $k$-1-($j$+$s$) least significant bits, and those of 2 add, with
saturation at the chosen maximum precision $k$.
\\[1ex]
The input precision $k$ must be taken equal to the desired output precision.
For instance, for an 8 x 8 bit multiplier with 16-bit output, odd input operands are
encoded as index $i$ in a 16-bit power $3^i$.
~~Addition is difficult in this code, so application is suggested for environments
restricted to multiplication mod $2^k$.

\subsection{Signed magnitude binary code over bases 2 and 3}

The proposed new number code is a signed magnitude code, well suited
for multiplication, and it uses two bases, namely 2 and 3. As shown, each k-digit
binary coded residue $n$ (mod $2^k$) is the product of a power $2^j$ of 2
~($j \leq k$), called the {\it even part} of $n$, and an odd residue called the
{\it odd part} of $n$, as shown the binary residue of a signed power
$\pm 3^i$ of 3 with $i < 2^{k-2}$.
\\[1ex]
Exponent pair ($i,j$) and sign $s$ uniquely encode each nonzero residue from
$-(2^k-1)$ to $2^k-1$, while the zero number 0 requires $j=k$, which can be
considered as an extra zero-bit $z$. \\[1ex]
To represent all $k$-bit binary numbers $n$ (integers), of which there are
$2^k$, a 4-component code $n=[z,s,t,u]$ is proposed, with the next interpretation:
\\[1ex] \hspace*{1cm}
$z$ : one {\bf zero} bit, with $z=0$ if $n=0$ and $z=1$ if $n \neq 0$.\\
\hspace*{1cm} $s$ : one {\bf sign} bit, with $s=0$ if $n>0$ and $s=1$ if $n<0$.
\\ \hspace*{1cm} $t$ : $k-2$ bits for the exponent $t$ of {\bf odd part} $3^t$.
\\ \hspace*{1cm}
$u$ : ~~$e$ bits for the exponent $u$ of {\bf even part} $2^u ~(u<k \leq 2^e$).
\\[1ex]
Extra {\bf overflow} bit ~$v=1$ iff $u_a+u_b \geq k$ : in case a product
$a.b$ exceeds $2^{k-1}$ in magnitude.
\\[1ex]
The code of the product of two such coded numbers $a=[z_a,s_a,t_a,u_a]$ and
$b=[z_b,s_b,t_b,u_b]$ is obtained by adding in binary code, by known means,
the odd and even code parts $t$ and $u$ respectively, and adding the signs
$s_a+s_b$ mod 2 (XOR), while multiplying the two zero bits $z_a.z_b$ (AND).
The overflow result bit $v=1$ iff the even part overflows: $u_a+u_b \geq k$.
\\[1ex]
Using for instance the known 'ripple-carry' way of binary addition hardware with
a full-adder cell FA per bit position, the schematic diagram is as follows,
where $t,t_a,t_b,u,u_a,u_b$ consist of 3 bits (of weights $2^0,2^1,2^2$),
and the optional overflow bit $v=u[2]*(u[1]+u[0])$ ~so iff $u \geq 5$:
\\[10ex]
\begin{picture}(50,50)
\setlength{\unitlength}{2mm}
\put( 5,14){$z_b$}              \put( 8,14){$z_a$}
\put( 6,13){\vector(0,-1){4}}   \put( 8,13){\vector(0,-1){4}}
\put( 5,5){\framebox(4,4){and}} \put( 7, 5){\vector(0,-1){3}}
\put( 7,1){$z$}
\put(13,14){$s_b$}              \put(16,14){$s_a$}
\put(14,13){\vector(0,-1){4}}   \put(16,13){\vector(0,-1){4}}
\put(13,5){\framebox(4,4){xor}} \put(15, 5){\vector(0,-1){3}}
\put(15,1){$s$}
\put(25,5){\framebox(4,4){$2^2$}} \put(27, 5){\vector(0,-1){3}}
\put(26,13){\vector(0,-1){4}}     \put(28,12){\vector(0,-1){3}}
\put(33,7){\vector(-1,0){4}}
\put(33,5){\framebox(4,4){$2^1$}} \put(35, 5){\vector(0,-1){3}}
\put(34,13){\vector(0,-1){4}}     \put(36,12){\vector(0,-1){3}}
\put(41,7){\vector(-1,0){4}}
\put(41,5){\framebox(4,4){$2^0$}} \put(43, 5){\vector(0,-1){3}}
\put(42,13){\vector(0,-1){4}}     \put(44,12){\vector(0,-1){3}}
\put(46,13){\line(-1,0){20}}  \put(47,13){$t_b$}
\put(46,12){\line(-1,0){18}}  \put(47,11){$t_a$}
\put(27,2){\vector(1,0){19}}  \put(47,1){$t$}
\put(55,5){\framebox(4,4){$2^2$}} \put(57, 5){\vector(0,-1){3}}
\put(56,13){\vector(0,-1){4}}     \put(58,12){\vector(0,-1){3}}
\put(63,7){\vector(-1,0){4}}
\put(63,5){\framebox(4,4){$2^1$}} \put(65, 5){\vector(0,-1){3}}
\put(64,13){\vector(0,-1){4}}     \put(66,12){\vector(0,-1){3}}
\put(71,7){\vector(-1,0){4}}
\put(71,5){\framebox(4,4){$2^0$}} \put(73, 5){\vector(0,-1){3}}
\put(72,13){\vector(0,-1){4}}     \put(74,12){\vector(0,-1){3}}
\put(76,13){\line(-1,0){20}}  \put(77,13){$u_b$}
\put(76,12){\line(-1,0){18}}  \put(77,11){$u_a$}
\put(57,2){\vector(1,0){19}}  \put(77,1){$u$}
\put(0,-3){{\bf Fig.1}:
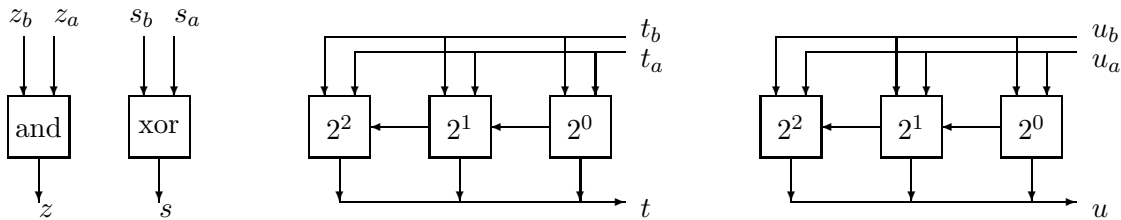
~Example multiplier mod $32=2^5$, with code $\pm~3^t.2^u~(t<2^3,~u \leq 5$) }
\end{picture}
\\[10ex]
{\bf Reference}: ~1. K.Hwang: {\it Computer Arithmetic}, J.Wiley \& Sons, NY 1979.

\end{document}